\documentclass[letterpaper]{IEEEtran}
\usepackage{latexsym,,amssymb,amsmath,graphicx,epsf,cite,bbm,float}
\usepackage{ifpdf}
\usepackage{epstopdf}

\usepackage{algorithm,algorithmic}
\usepackage{amsmath,amssymb,bm}
\usepackage{amsfonts,dsfont,color,bbm,subcaption}

\usepackage{amsmath,amsthm}
\usepackage{hyperref}

\newtheorem{theorem}{Theorem}

\newtheorem{lemma}[]{Lemma}

\newtheorem{remark}[]{Remark}

\newtheorem{definition}[]{Definition}

\newcommand{\be}{\begin{equation}}
\newcommand{\ee}{\end{equation}}
\newcommand{\bea}{\begin{eqnarray}}
\newcommand{\eea}{\end{eqnarray}}

\newcommand{\MB}{\left[\begin{array}}
\newcommand{\ME}{\end{array}\right]}

\newcommand{\ei}{\end{itemize}}
\newcommand{\bi}{\begin{itemize}}

\renewcommand{\Re}{\mathbb{R}}

\usepackage{amsmath,amsthm}

\title{$1D$ to $nD$: A Meta Algorithm for Multivariate Global Optimization via Univariate Optimizers}
\author{\IEEEauthorblockN{Kaan Gokcesu}, \IEEEauthorblockN{Hakan Gokcesu} }
\begin{document}
\maketitle

\begin{abstract}
	We propose a meta algorithm that can solve a multivariate global optimization problem using univariate global optimizers. Although the univariate global optimization does not receive much attention compared to the multivariate case, which is more emphasized in academia and industry; we show that it is still relevant and can be directly used to solve problems of multivariate optimization. We also provide the corresponding regret bounds in terms of the time horizon $T$ and the average regret of the univariate optimizer, when it is robust against bounded noises with robust regret guarantees.   
\end{abstract}

\section{Introduction}

Global optimization aims to minimize the value of a loss function with minimal number of evaluations, which is paramount in a myriad of applications \cite{cesa_book, poor_book}. In these kinds of problem scenarios, the performance of a parameter set requires to be evaluated numerically or with cross-validations, which may have a high computational cost. Furthermore, since the search space needs to be sequentially explored, the number of samples to be evaluated needs to be small because of various operational constraints. The evaluation samples may be hard to determine especially if the function does not possess common desirable properties such as convexity or linearity. 
	
	This sequential optimization problem of unknown and possibly non-convex functions are referred to as global optimization \cite{pinter1991global}. It has also been dubbed by the names of derivative-free optimization \cite{rios2013derivative} (since the optimization is done based on solely the function evaluations and the derivatives are inconsequential) and black-box optimization \cite{jones1998efficient}. Over the past years, the global optimization problem has gathered significant attention with various algorithms being proposed in distinct fields of research. It has been studied especially in the fields of non-convex optimization \cite{jain2017non, hansen1991number, basso1982iterative}, Bayesian optimization \cite{brochu2010tutorial}, convex optimization \cite{boyd2004convex,nesterov2003introductory,bubeck2015convex}, bandit optimization \cite{munos2014bandits}, stochastic optimization \cite{shalev2012online, spall2005introduction}; because of its practical applications in distribution estimation \cite{willems,gokcesu2017online,coding2}, multi agent systems \cite{gokcesu2020generalized,vural2019minimax,cesa-bianchi,gokcesu2018bandit}, control theory \cite{tnnls3}, robust learning \cite{gokcesu2021generalized}, signal processing \cite{ozkan}, prediction \cite{gokcesu2020recursive,singer}, game theory \cite{tnnls1}, decision theory \cite{gokcesu2021optimally,tnnls4} and anomaly detection \cite{gokcesu2019outlier}.
	
	Global optimization studies the problem of optimizing a function $f(\cdot)$ with minimal number of evaluations. Generally, only the evaluations of the sampled points are revealed and the information about $f(\cdot)$ or its derivatives are unavailable. When a point $x$ is queried, only its value $f(x)$ is revealed. 
	There exists various heuristic algorithms in literature that deals with this problem such as model-based methods, genetic algorithms and Bayesian optimization. However, most popular approach is still the regularity-based methods since in many applications the system has some inherent regularity with respect to its input-output pair, i.e., $f(\cdot)$ satisfies some regularity condition. 
	
	Even though some works such as \cite{bartlett2019simple} and \cite{grill2015black} use a smoothness regularity in regards to the hierarchical partitioning; traditionally, Lipschitz continuity or smoothness are more common in the literature. 
	Lipschitz regularity was first studied in the works of \cite{piyavskii1972algorithm, shubert1972sequential}, where the algorithm has been dubbed the Piyavskii-Shubert algorithm. A lot of research has been done with this algorithm at its base such as \cite{jacobsen1978global,basso1982iterative,mayne1984outer,mladineo1986algorithm,shen1987interval,horst1987convergence,hansen1991number}, which study special aspects of its application with examples involving functions that satisfy some Lipschitz condition and propose alternative formulations.
	
	Breiman and Cutler \cite{breiman1993deterministic} proposed a multivariate extension that utilizes the Taylor expansion of $f(\cdot)$ at its core. 
	Baritompa and Cutler \cite{baritompa1994accelerations} propose an alternative acceleration of the Breiman and Cutler’s method. 
	Hansen and Jaumard \cite{hansen1995lipschitz} summarize and discuss the algorithms in literature. 
	Sergeyev \cite{sergeyev1998global} proposes a variant that utilizes smooth auxiliary functions. 
	Ellaia et al. \cite{ellaia2012modified} suggest a variant that maximizes a univariate differentiable function.
	Brent \cite{brent2013algorithms} proposes a variant where the function $f(\cdot)$ is required to be defined on a compact interval with a bounded second derivative. The works in \cite{gokcesu2021regret,gokcesu2022low} study the problem with more generalized Lipschitz regularities.
	Horst and Tuy \cite{horst2013global} provide a discussion about the role of deterministic algorithms in global optimization.
	
	
	The study of the number of iterations of Piyavskii-Shubert algorithm was initiated by Danilin \cite{danilin1971estimation}. For its simple regret analysis, a crude regret bound of the form $\tilde{r}_T = O(T^{-1})$ can be obtained from \cite{mladineo1986algorithm} when the function $f(\cdot)$ is Lipschitz continuous. 
	For univariate functions, the work by Hansen et al. \cite{hansen1991number} derives a bound on the sample complexity for a Piyavskii–Shubert variant, which stops automatically upon returning an $\epsilon$-optimizer. They showed that the number of evaluations required is at most proportional with $\int_0^1(f(x^*)-f(x)+\epsilon)^{-1}dx$, which improves upon the results of \cite{danilin1971estimation}. 
	The work by Ellaia et al. \cite{ellaia2012modified} improves upon the previous results of \cite{danilin1971estimation,hansen1991number}.
	The work of Malherbe and Vayatis \cite{malherbe2017global} studies a variant of the Piyavskii–Shubert algorithm and obtain regret upper bounds under strong assumptions. The work of Bouttier et al. \cite{bouttier2020regret} studies the regret bounds of Piyavskii-Shubert algorithm under noisy evaluations.
	Instead of the weaker simple regret, the work in \cite{gokcesu2021regret} provides cumulative regret bounds for variants of Piyavskii-Shubert algorithm. 
	The work in \cite{gokcesu2022low} provides a more applicable and computationally superior alternative with similar regret guarantees.
	
	While the problem of global optimization is most meaningful in multivariate settings, many works deal with univariate objectives because of the simplicity of the problem setting, low computational complexity and better performance guarantees. To this end, we propose a meta algorithm that can utilize univariate optimizers to solve multivariate objectives.
	
	\section{Preliminaries and the Meta Algorithm}\label{sec:problem}
	In this section, we provide the formal problem definition. In the multivariate global optimization problem, we want to optimize a function $f(\cdot)$ that maps from $d$-dimensional unit cube $[0,1]^d$ to the real line $\Re$, i.e.,
	\begin{align}
		f(\cdot) : [0,1]^d\rightarrow\Re,\label{eq:f}
	\end{align}
	
	However, it is not straightforward to optimize any arbitrary function $f(\cdot)$. To this end, we define a regularity measure. Instead of the restrictive Lipschitz continuity or smoothness \cite{gokcesu2021regret}; we define a weaker, more general regularity condition.
	\begin{definition}\label{def:condition}
		Let the function $f(\cdot)$ that we want to optimize satisfy the following condition:
		\begin{align*}
			|f(x)-f(y)|\leq d({x,y}),
		\end{align*}
		where $d(\cdot)$ is a known norm induced metric.
	\end{definition}

	We optimize the function $f(\cdot)$ iteratively by selecting a query point $x_t$ at each time $t$ and receive its evaluation $f(x_t)$. Then, we select the next query point based on the past queries and their evaluations. Hence,
	\begin{align}
		x_{t+1}=\Gamma(x_1,x_2,\ldots,x_{t},f(x_1),f(x_2),\ldots,f(x_{t})),
	\end{align}
	where $\Gamma(\cdot)$ is some function. One such algorithm is the famous Piyavskii–Shubert algorithm (and its variants) \cite{gokcesu2021regret}.
	
	We approach this problem from a loss minimization perspective (as in line with the computational learning theory). We consider the objective function $f(\cdot)$ as a loss function to be minimized by producing query points (predictions) $x_t$ from the compact subset $\Theta$ at each point in time $t$. We define the performance of the predictions $x_t$ for a time horizon $T$ by the cumulative loss incurred up to $T$ instead of the best loss so far at time $T$; i.e., instead of the loss of the best prediction up to time $T$:
	$\tilde{l}_t \triangleq \min_{t\in\{1,\ldots,T\}}f(x_t),$
	we use the average loss up to time $T$:
	$l_T\triangleq\frac{1}{T}\sum_{t=1}^{T}f(x_t).$
	Let $f_*$ be a global minimum of $f(\cdot)$, i.e.,
	$f_* = \min_{x\in\Theta} f(x).\label{eq:min}
	$
	As in line with learning theory, we use the notion of regret to evaluate the performance of our algorithm \cite{gokcesu2021regret}. Hence, instead of the simple regret at time $T$, we analyze the average regret up to time $T$:
	\begin{align}
		r_T\triangleq\frac{1}{T}\sum_{t=1}^{T}f(x_t)-f_*\label{eq:regret}.
	\end{align}

Let us have a robust $1D$ optimizer that can produce average regret guarantees when the evaluations are noisy. Given a univariate objective $h(x)$, let us observe a noisy version $\tilde{h}(x)=h(x)+\epsilon(x)$, where the noise is nonnegative and bounded, i.e.,  $0\leq \epsilon(x)\leq \epsilon$ for some known $\epsilon$, possibly after a bias translation.

Given the noisy function evaluations $\{\tilde{h}(x_{\tau})\}_{\tau=1}^{t-1}$, noise bound $\epsilon$ and time horizon $T$; let the query at time $t$ of the $1D$ optimizer be the output of the following function 
\begin{align}
	x_t=\Gamma_{T,\epsilon}^t(x_{1},x_{2},\ldots,x_{t-1}; \tilde{h}(x_{1}),\tilde{h}(x_{2}),\ldots,\tilde{h}(x_{t-1})).
\end{align}
To solve a $d$-dimensional objective $f(\cdot)$, we use a meta algorithm that takes the $1D$ optimizer together with parameters $\{\epsilon_i\}_{i=1}^d$ and $\{T_i\}_{i=1}^d$ as inputs. The meta algorithm will query $T_1$ points in the first dimension and for each individual query at dimension $i$, it will have $T_{i+1}$ queries at dimension $i+1$.

The $d$-dimensional query at time $t$ is given by some $\{x^i_{\tau_{i}}\}_{i=1}^d$, where $\tau_i$ are such that
\begin{align}
	t=1+\sum_{i=1}^d(\tau_i-1)\prod_{j=i+1}^{d}T_j.
\end{align}
Using the $1D$ optimizer, the queries at each dimension are chosen with
\begin{align}
	x^i_{\tau_{i}}=\Gamma^{\tau_i}_{T_i,\epsilon_i}(\{x^i_{j}\}_{j=1}^{\tau_i-1};h^i(\{x^i_{j}\})_{j=1}^{\tau_i-1}),
\end{align}
where the functions $h^i(x^i_j)$ is the minimum evaluation made when the $i^{th}$ dimension query is $x^{i}_j$. We point out that the meta-algorithm solves the minimum objective in each dimension, hence, the regularity at \autoref{def:condition} is preserved. However, instead of the minimum possible evaluations, we only have access to the minimum evaluation we make so far, which why the nonnegative noise comes in to play. 

A pseudo-code is provided in \autoref{alg}.
	
\begin{algorithm}[!t]
	\caption{Meta Algorithm}\label{alg}
	{\begin{algorithmic}[1]
		\STATE Inputs $\{\epsilon_i\}_{i=1}^d$, $\{T_i\}_{i=1}^d$.
		\STATE $d$-dimensional objective $f(\{x^i\})$
		\FOR{$i=1$ \TO $d$}
			\STATE $\tau_i=1$
		\ENDFOR
		\FOR{$i=1$ \TO $d$}
		\STATE Set initial query $x^i_{\tau_i}$
		\ENDFOR
		\STATE Evaluate $f_1=f(\{x_{\tau_i}^i\}_{i=1}^d)$
		\FOR{$i=1$ \TO $d$}
		\STATE $h^i(x^i_{\tau_{i}})=f_1$
		\ENDFOR
		\FOR{$t=2$ \TO $\prod_{i=1}^d T_i$}
			\STATE $\tau_d\leftarrow \tau_d+1$
			\FOR{$i=d$ \TO $2$}
				\IF{$\tau_i>T_i$}
				\STATE	$\tau_i\leftarrow 1$
				\STATE $\tau_{i-1}\leftarrow\tau_{i-1}+1$
				\ENDIF	
			\ENDFOR
			\STATE Set $I\leftarrow d$
			\WHILE{$\tau_I=1$}
				\STATE Set initial query $x^I_{\tau_I}$
				\STATE $I\leftarrow I-1$
			\ENDWHILE
			\STATE Set $ x^I_{\tau_I}\leftarrow\Gamma_{T_I,\epsilon_I}^{\tau_I}(\{x^I_\tau\}_{\tau=1}^{\tau_I-1};\{h^I(x^I_\tau)\}_{\tau=1}^{\tau_I-1})$
			\STATE Evaluate $f_t=f(\{x_{\tau_i}^i\}_{i=1}^d)$
			\FOR{$i=d$ \TO $I$}
			\STATE Set $h^i(x^i_{\tau_i})\leftarrow f_t$
			\ENDFOR
			\FOR{$i=I-1$ \TO $1$}
				\STATE $h^i(x^i_{\tau_i})\leftarrow \min(f_t,h^{i}(x^{i}_{\tau_{i}}))$
			\ENDFOR
		\ENDFOR
	\end{algorithmic}}
\end{algorithm}	
		
\section{Regret Analyses}
Let the meta algorithm make $T_x=T$ number of predictions $x_{i,n}=[y_{i,n}\: z_n]$, where $i\in\{1,\ldots,T_y\}$ and $n\in \{1,\ldots,T_z\}$ for some $T_x=T=T_yT_z$.
The average regret will be given by
\begin{align}
	r_{T_x}^{(d_x)}=\frac{1}{T_x}\sum_{i,n}f(x_{i,n})-f_*,
\end{align}
where $d_x=d$ is the dimension of $x$ (i.e., $d_x=d_y+d_z$) and $f_*=\min_x f(x)$ is the optimal evaluation. Suppose instead of $f(x_{i,n})$, we observe a noisy $\tilde{f}(x_{i,n})$ such that
\begin{align}
	\tilde{f}(x_{i,n})=f(x_{i,n})+\epsilon(x_{i,n}),
\end{align}
where $0\leq\epsilon(x_{i,n})\leq \epsilon_x$ for some known $\epsilon_x\in[0,\infty]$. Let us define the following quantities when $\tilde{f}(\cdot)$ is observed instead of $f(\cdot)$.
\begin{itemize}
	\item We define the average regret bound as
	\begin{align}
		r_{T_x}^{(d_x)}(\epsilon_x)\geq \frac{1}{T_x}\sum_{i,n}f(x_{i,n})-f_*.
	\end{align}
	\item We also define the average pseudo-regret bound as
	\begin{align}
		\tilde{r}_{T_x}^{(d_x)}(\epsilon_x)\geq\frac{1}{T_x}\sum_{i,n}\tilde{f}(x_{i,n})-f_*.
	\end{align}
\end{itemize}
We assume that the average bounds are nonincreasing with $T$. We have the following recursive relation.
\begin{lemma}\label{thm:regretRec}
For the predictions $x_{i,n}=[y_{i,n}\: z_n]$ by the meta-algorithm, we have the following average regret relation:
	\begin{align*}
		r_{T_x}^{(d_x)}(\epsilon_x)\leq\:&r_{T_y}^{(d_y)}(\epsilon_x)+r^{(d_z)}_{T_z}(\epsilon_x+\epsilon_y),
	\end{align*}
	where $\epsilon_y=\max_{n}(\min_i f([y_{i,n}\:z_n])-\min_y f([y\:z_n]))$.
\begin{proof}
	We have
	\begin{align*}
		r_{T_x}^{(d_x)}(\epsilon_x)=&\frac{1}{T_x}\sum_{i,n}f(x_{i,n})-f_*,
		\\=&\frac{1}{T_x}\sum_{i,n} f\left([y_{i,n}\: z_n]\right)-f_*
		\\=&\frac{1}{T_x}\sum_{i,n} f\left([y_{i,n}\: z_n]\right)-\frac{1}{T_z}\sum_{n} \min_y f\left([y\: z_n]\right)
		\\&+\frac{1}{T_z}\sum_{n} \min_y f\left([y\: z_n]\right)-f_*,
		\\=&\frac{1}{T_z}\sum_{n}\left[\frac{1}{T_y} \sum_i f\left([y_{i,n}\: z_n]\right)- \min_y f\left([y\: z_n]\right)\right]
		\\&+\frac{1}{T_z}\sum_{n} \min_y f\left([y\: z_n]\right)-f_*,
		\\\leq&\max_n \left\{\frac{1}{T_y} \sum_i f\left([y_{i,n}\: z_n]\right)- \min_y f\left([y\: z_n]\right)\right\}
		\\&+\frac{1}{T_z}\sum_{n} \min_y f\left([y\: z_n]\right)-f_*,
		\\\leq&r_{T_y}^{(d_y)}(\epsilon_x)+r^{(d_z)}_{T_z}(\epsilon_x + \epsilon_y),
	\end{align*}
	where $\epsilon_y=\max_{n}(\min_i f([y_{i,n}\:z_n])-\min_y f([y\:z_n]))$.
\end{proof}
\end{lemma}

\subsection{Strongly robust}
Suppose our $1D$ optimizer is strongly robust in the sense that given the noisy evaluations
\begin{align}
	\tilde{h}(x_{t})=h(x_t)+\epsilon(x_t)
\end{align}
instead of $h(x_t)$, it can still guarantee a regret bound. Let it achieve an average regret bound of
\begin{align}
	{r}^{(1)}_T(\epsilon)\leq r^{(1)}_T(0)+\alpha_T\epsilon,\label{eq:stronglyrobust}
\end{align}
where $\epsilon\geq \epsilon(x_t)\geq 0$ for some non-increasing $\alpha_T\geq 0$.

\begin{lemma}\label{thm:recRobust}
	We have
	\begin{align}
		r_{T_x}^{(d_x)}(0)\leq\:&(1+\alpha_{T_z})r_{T_y}^{(d_x-1)}(0)+r^{(1)}_{T_z}(0),
	\end{align}
	\begin{proof}
		Using \autoref{thm:regretRec} with the fact that $\epsilon_y$
		is upper-bounded by the average regret, we have  
		\begin{align}
			r_{T_x}^{(d_x)}(\epsilon_x)\leq\:&r_{T_y}^{(d_y)}(\epsilon_x)+r^{(d_z)}_{T_z}(\epsilon_x+\epsilon_y),
			\\\leq&r_{T_y}^{(d_y)}(\epsilon_x)+r^{(d_z)}_{T_z}(\epsilon_x+r_{T_y}^{(d_y)}(\epsilon_x)).
		\end{align}
	Setting $\epsilon_x=0$ gives
	\begin{align}
		r_{T_x}^{(d_x)}(0)\leq\:&r_{T_y}^{(d_y)}(0)+r^{(d_z)}_{T_z}(r_{T_y}^{(d_y)}(0)).
	\end{align}
	Let $d_z=1$ and $d_y=d_x-1$, which gives
	\begin{align}
		r_{T_x}^{(d_x)}(0)\leq\:&r_{T_y}^{(d_x-1)}(0)+r^{(1)}_{T_z}(r_{T_y}^{(d_x-1)}(0)).
	\end{align}
	Using \eqref{eq:stronglyrobust} provides
	\begin{align}
		r_{T_x}^{(d_x)}(0)\leq\:&r_{T_y}^{(d_x-1)}(0)+\alpha_{T_z}r_{T_y}^{(d_x-1)}(0)+r^{(1)}_{T_z}(0)
		\\\leq&(1+\alpha_{T_z})r_{T_y}^{(d_x-1)}(0)+r^{(1)}_{T_z}(0),
	\end{align}
	which concludes the proof.
	\end{proof}
\end{lemma}

\begin{theorem}\label{thm:regretStrong}
	Let the meta-algorithm produce its $T=\prod_{i=1}^d T_i$ number of $d$-dimensional predictions, where $T_1\leq T_2\leq \ldots\leq T_d$. We have the following average regret
	\begin{align}
		r_{\prod_{i=1}^d T_i}^{(d)}(0)\leq&d(1+\alpha_{T_1})^{d-1}r^{(1)}_{T_{1}}(0).
	\end{align}
	\begin{proof}
		Using \autoref{thm:recRobust}, we have
	\begin{align}
		r_{\prod_{i=j}^d T_i}^{(d-j+1)}(0)\leq&(1+\alpha_{T_j})r_{\prod_{i=j+1}^d T_i}^{(d-j)}(0)+r^{(1)}_{T_j}(0).
	\end{align}	
	Hence, the recursion gives
	\begin{align}
		r_{\prod_{i=1}^d T_i}^{(d)}(0)\leq&r^{(1)}_{T_1}(0)+\sum_{i=1}^{d-1}r^{(1)}_{T_{i+1}}(0)\prod_{j=1}^{i}(1+\alpha_{T_j}),
		\\\leq&r^{(1)}_{T_1}(0)+r^{(1)}_{T_{1}}(0)\sum_{i=1}^{d-1}\prod_{j=1}^{i}(1+\alpha_{T_j}),
		\\\leq&r^{(1)}_{T_1}(0)+r^{(1)}_{T_{1}}(0)\sum_{i=1}^{d-1}(1+\alpha_{T_1})^i,
		\\\leq&r^{(1)}_{T_{1}}(0)\sum_{i=0}^{d-1}(1+\alpha_{T_1})^i,
		\\\leq&d(1+\alpha_{T_1})^{d-1}r^{(1)}_{T_{1}}(0),
	\end{align}
	which concludes the proof.	 
	\end{proof}
\end{theorem}

\subsection{Weakly robust}
Suppose our $1D$ optimizer is weakly robust in the sense that given the noisy evaluations
\begin{align}
	\tilde{h}(x_{t})=h(x_t)+\epsilon(x_t)
\end{align}
instead of $h(x_t)$, it can still guarantee a regret bound. Let it achieve an average pseudo-regret bound of
\begin{align}
	\tilde{r}^{(1)}_T(\epsilon)\leq r^{(1)}_T(0)+\beta_T\epsilon,\label{eq:weaklyrobust}
\end{align}
where $\epsilon\geq \epsilon(x_t)\geq 0$ for some non-increasing $\beta_T\geq 1$. We call this scenario weakly robust, since the strongly robust guarantee directly implies a weakly robust guarantee with $\beta_T=\alpha_T+1$.

\begin{lemma}\label{thm:recWRobust}
	We have
	\begin{align}
	r_{T_x}^{(d_x)}(\epsilon_x)\leq\:&\tilde{r}_{T_y}^{(1)}(\epsilon_x)+r^{(d_x-1)}_{T_z}(\tilde{r}_{T_y}^{(1)}(\epsilon_x)),
\end{align}
	\begin{proof}
		Using the fact that average $\epsilon_x(\cdot)+\epsilon_y(\cdot)$ is upper-bounded by the average pseudo-regret, we have, instead of \autoref{thm:regretRec},  
		\begin{align}
			r_{T_x}^{(d_x)}(\epsilon_x)
			\leq&r_{T_y}^{(d_y)}(\epsilon_x)+r^{(d_z)}_{T_z}(\tilde{r}_{T_y}^{(d_y)}(\epsilon_x)).
		\end{align}
		Let $d_z=d_x-1$ and $d_y=1$, which gives
		\begin{align}
			r_{T_x}^{(d_x)}(\epsilon_x)\leq\:&r_{T_y}^{(1)}(\epsilon_x)+r^{(d_x-1)}_{T_z}(\tilde{r}_{T_y}^{(1)}(\epsilon_x)),
			\\\leq\:&\tilde{r}_{T_y}^{(1)}(\epsilon_x)+r^{(d_x-1)}_{T_z}(\tilde{r}_{T_y}^{(1)}(\epsilon_x)),
		\end{align}
		which concludes the proof.
	\end{proof}
\end{lemma}

\begin{theorem}\label{thm:regretWeak}
	Let the meta-algorithm produce its $T=\prod_{i=1}^d T_i$ number of $d$-dimensional predictions, where $T_1\leq T_2\leq \ldots\leq T_d$. We have the following average regret
	\begin{align}
		r_{\prod_{i=1}^d T_i}^{(d)}(0)\leq&0.5(d+1)d\beta_{T_{1}}^{d-1}r^{(1)}_{T_1}(0),
	\end{align}
	\begin{proof}
		Using \autoref{thm:recWRobust}, we have
		\begin{align}
			r_{\prod_{i=j}^d T_i}^{(d-j+1)}(\epsilon)\leq&\tilde{r}^{(1)}_{T_j}(\epsilon)+r_{\prod_{i=j+1}^d}^{(d-j)}(\tilde{r}^{(1)}_{T_j}(\epsilon)).
		\end{align}	
		Hence, the recursion gives
		\begin{align}
			r^{(d)}_{\prod_{i=1}^d T_i}(0)\leq \sum_{j=1}^d \tilde{r}^{(1)}_{T_j}\circ\tilde{r}^{(1)}_{T_{j-1}}\circ\ldots\circ\tilde{r}^{(1)}_{T_2}\circ\tilde{r}^{(1)}_{T_1}(0).
		\end{align}
		Using \eqref{eq:weaklyrobust}, we have
		\begin{align}
			r^{(d)}_{\prod_{i=1}^d T_i}(0)\leq& \sum_{j=1}^d \sum_{k=1}^{j}r^{(1)}_{T_k}(0)\prod_{l=k+1}^{j}\beta_{T_l}
			\\\leq&\sum_{j=1}^d \sum_{k=1}^{j}r^{(1)}_{T_k}(0)\beta_{T_{k+1}}^{j-k}
			\\\leq&\sum_{j=1}^d \sum_{k=1}^{j}r^{(1)}_{T_k}(0)\beta_{T_{1}}^{j-k}
			\\\leq&\sum_{j=1}^dr^{(1)}_{T_1}(0) \sum_{k=0}^{j-1}\beta_{T_{1}}^{k}
			\\\leq&r^{(1)}_{T_1}(0)\sum_{j=0}^{d-1} \sum_{k=0}^{j}\beta_{T_{1}}^{k}
			\\\leq&\frac{d(d+1)}{2}\beta_{T_{1}}^{d-1}r^{(1)}_{T_1}(0),
		\end{align}
		which concludes the proof.	 
	\end{proof}
\end{theorem}

\subsection{Arbitrary $T$}
$T$ may not necessarily be factorisable. To this end, for any arbitrary $T$, the meta algorithm selections are as follows:

\begin{remark}\label{rem:selection}
	Given input $T$, the meta algorithm chooses $T_i\in\{\lfloor T^{\frac{1}{d}}\rfloor   , \lceil T^{\frac{1}{d}}\rceil \},$ $\forall i$ such that $T_i$ are in nonincreasing order, i.e., $T_1\leq T_2\leq\ldots\leq T_d$, and $T\leq \prod_{i=1}^d T_i \leq \frac{\lceil T^{\frac{1}{d}}\rceil}{\lfloor T^{\frac{1}{d}}\rfloor}T$.	
\end{remark}

\begin{remark}\label{rem:regret}
We observe that in both \autoref{thm:regretStrong} and \autoref{thm:regretWeak}, the regret is in the following form
\begin{align}
	r^{(d)}_{\prod_{i=1}^dT_i}(0)\leq \mathcal{F}(d,T_1)r^{(1)}_{T_1}(0),
\end{align}
for some function $\mathcal{F}(d,T_1)$, which is either $d(1+\alpha_{T_1})^{d-1}$ or $0.5(d+1)d\beta_{T_1}^{d-1}$ depending on the robustness level of the optimizer. 
\end{remark}

\begin{lemma}\label{thm:cumReg}
	With the selections in \autoref{rem:selection}, we have the following cumulative regret
	\begin{align}
		R^{(d)}_T(0)\leq&2T\mathcal{F}(d,\lfloor T^{\frac{1}{d}}\rfloor)r^{(1)}_{\lfloor T^{\frac{1}{d}}\rfloor}(0).
	\end{align}
	\begin{proof}
		Since the instantaneous regret is always non-negative, we have the following for cumulative regrets
		\begin{align}
			R^{(d)}_T(0)\leq R^{(d)}_{\prod_{i=1}^dT_i}(0),
		\end{align} 
		when $T_i$ are selected with \autoref{rem:selection}. From \autoref{rem:regret}, we have
		\begin{align}
			R^{(d)}_T(0)\leq& R^{(d)}_{\prod_{i=1}^dT_i}(0),
			\\\leq&\mathcal{F}(d,T_1)r^{(1)}_{T_1}(0)\prod_{i=1}^{d}T_i
			\\\leq&2T\mathcal{F}(d,T_1)r^{(1)}_{T_1}(0)
			\\\leq&2T\mathcal{F}(d,\lfloor T^{\frac{1}{d}}\rfloor)r^{(1)}_{\lfloor T^{\frac{1}{d}}\rfloor}(0),
		\end{align} 
		which concludes the proof.
	\end{proof}
\end{lemma}
For unknown time horizon, we utilize the doubling trick and reset the meta algorithm in each epoch. Let $T=\sum_{n=0}^{N-1}2^n+C$, where $1\leq C\leq 2^N$. Hence, $2^N\leq T\leq 2^{N+1}\leq 2T$. 
\begin{theorem}
	For unknown horizon $T$ and the selections in \autoref{rem:selection}, we have the following average regret
	\begin{align}
		\bar{r}^{(d)}_T(0)\leq&2\log_2(2T)\mathcal{F}(d,\lfloor T^{\frac{1}{d}}\rfloor)r^{(1)}_{\lfloor T^{\frac{1}{d}}\rfloor}(0).
	\end{align}
	\begin{proof}
		From \autoref{thm:cumReg}, we will have
		\begin{align}
			\bar{R}^{(d)}_T(0)\leq&\sum_{n=0}^{N-1}R^{(d)}_{2^n}(0)+R^{(d)}_C(0)
			\\\leq&(N+1)R^{(d)}_{2^N}(0)
			\\\leq&(N+1)R^{(d)}_{T}(0)
			\\\leq&2\log_2(2T)T\mathcal{F}(d,\lfloor T^{\frac{1}{d}}\rfloor)r^{(1)}_{\lfloor T^{\frac{1}{d}}\rfloor}(0)
		\end{align}
		Hence, the average regret is bounded as
		\begin{align}
			\bar{r}^{(d)}_T(0)\leq&2\log_2(2T)\mathcal{F}(d,\lfloor T^{\frac{1}{d}}\rfloor)r^{(1)}_{\lfloor T^{\frac{1}{d}}\rfloor}(0),
		\end{align}
		which concludes the proof.
	\end{proof}
\end{theorem}

\bibliographystyle{IEEEtran}
\bibliography{double_bib}

\begin{thebibliography}{10}
\providecommand{\url}[1]{#1}
\csname url@samestyle\endcsname
\providecommand{\newblock}{\relax}
\providecommand{\bibinfo}[2]{#2}
\providecommand{\BIBentrySTDinterwordspacing}{\spaceskip=0pt\relax}
\providecommand{\BIBentryALTinterwordstretchfactor}{4}
\providecommand{\BIBentryALTinterwordspacing}{\spaceskip=\fontdimen2\font plus
\BIBentryALTinterwordstretchfactor\fontdimen3\font minus
  \fontdimen4\font\relax}
\providecommand{\BIBforeignlanguage}[2]{{%
\expandafter\ifx\csname l@#1\endcsname\relax
\typeout{** WARNING: IEEEtran.bst: No hyphenation pattern has been}%
\typeout{** loaded for the language `#1'. Using the pattern for}%
\typeout{** the default language instead.}%
\else
\language=\csname l@#1\endcsname
\fi
#2}}
\providecommand{\BIBdecl}{\relax}
\BIBdecl

\bibitem{cesa_book}
N.~Cesa-Bianchi and G.~Lugosi, \emph{Prediction, learning, and games}.\hskip
  1em plus 0.5em minus 0.4em\relax Cambridge university press, 2006.

\bibitem{poor_book}
H.~V. Poor, \emph{An Introduction to Signal Detection and Estimation}.\hskip
  1em plus 0.5em minus 0.4em\relax NJ: Springer, 1994.

\bibitem{pinter1991global}
J.~D. Pint{\'e}r, ``Global optimization in action,'' \emph{Scientific
  American}, vol. 264, pp. 54--63, 1991.

\bibitem{rios2013derivative}
L.~M. Rios and N.~V. Sahinidis, ``Derivative-free optimization: a review of
  algorithms and comparison of software implementations,'' \emph{Journal of
  Global Optimization}, vol.~56, no.~3, pp. 1247--1293, 2013.

\bibitem{jones1998efficient}
D.~R. Jones, M.~Schonlau, and W.~J. Welch, ``Efficient global optimization of
  expensive black-box functions,'' \emph{Journal of Global optimization},
  vol.~13, no.~4, pp. 455--492, 1998.

\bibitem{jain2017non}
P.~Jain and P.~Kar, ``Non-convex optimization for machine learning,''
  \emph{Foundations and Trends{\textregistered} in Machine Learning}, vol.~10,
  no. 3-4, pp. 142--336, 2017.

\bibitem{hansen1991number}
P.~Hansen, B.~Jaumard, and S.-H. Lu, ``On the number of iterations of
  piyavskii's global optimization algorithm,'' \emph{Mathematics of Operations
  Research}, vol.~16, no.~2, pp. 334--350, 1991.

\bibitem{basso1982iterative}
P.~Basso, ``Iterative methods for the localization of the global maximum,''
  \emph{SIAM Journal on Numerical Analysis}, vol.~19, no.~4, pp. 781--792,
  1982.

\bibitem{brochu2010tutorial}
E.~Brochu, V.~M. Cora, and N.~De~Freitas, ``A tutorial on bayesian optimization
  of expensive cost functions, with application to active user modeling and
  hierarchical reinforcement learning,'' \emph{arXiv preprint arXiv:1012.2599},
  2010.

\bibitem{boyd2004convex}
S.~Boyd, S.~P. Boyd, and L.~Vandenberghe, \emph{Convex optimization}.\hskip 1em
  plus 0.5em minus 0.4em\relax Cambridge university press, 2004.

\bibitem{nesterov2003introductory}
Y.~Nesterov, \emph{Introductory lectures on convex optimization: A basic
  course}.\hskip 1em plus 0.5em minus 0.4em\relax Springer Science \& Business
  Media, 2003, vol.~87.

\bibitem{bubeck2015convex}
S.~Bubeck, ``Convex optimization: Algorithms and complexity,''
  \emph{Foundations and Trends{\textregistered} in Machine Learning}, vol.~8,
  no. 3-4, pp. 231--357, 2015.

\bibitem{munos2014bandits}
R.~Munos, ``From bandits to monte-carlo tree search: The optimistic principle
  applied to optimization and planning,'' \emph{Foundations and Trends® in
  Machine Learning}, vol.~7, no.~1, pp. 1--129, 2014.

\bibitem{shalev2012online}
S.~Shalev-Shwartz \emph{et~al.}, ``Online learning and online convex
  optimization,'' \emph{Foundations and Trends{\textregistered} in Machine
  Learning}, vol.~4, no.~2, pp. 107--194, 2012.

\bibitem{spall2005introduction}
J.~C. Spall, \emph{Introduction to stochastic search and optimization:
  estimation, simulation, and control}.\hskip 1em plus 0.5em minus 0.4em\relax
  John Wiley \& Sons, 2005, vol.~65.

\bibitem{willems}
F.~M.~J. Willems, ``Coding for a binary independent
  piecewise-identically-distributed source.'' \emph{IEEE Transactions on
  Information Theory}, vol.~42, no.~6, pp. 2210--2217, 1996.

\bibitem{gokcesu2017online}
K.~Gokcesu and S.~S. Kozat, ``Online anomaly detection with minimax optimal
  density estimation in nonstationary environments,'' \emph{IEEE Transactions
  on Signal Processing}, vol.~66, no.~5, pp. 1213--1227, 2017.

\bibitem{coding2}
G.~I. Shamir and N.~Merhav, ``Low-complexity sequential lossless coding for
  piecewise-stationary memoryless sources,'' \emph{IEEE Transactions on
  Information Theory}, vol.~45, no.~5, pp. 1498--1519, Jul 1999.

\bibitem{gokcesu2020generalized}
K.~Gokcesu and H.~Gokcesu, ``A generalized online algorithm for translation and
  scale invariant prediction with expert advice,'' \emph{arXiv preprint
  arXiv:2009.04372}, 2020.

\bibitem{vural2019minimax}
N.~M. Vural, H.~Gokcesu, K.~Gokcesu, and S.~S. Kozat, ``Minimax optimal
  algorithms for adversarial bandit problem with multiple plays,'' \emph{IEEE
  Transactions on Signal Processing}, vol.~67, no.~16, pp. 4383--4398, 2019.

\bibitem{cesa-bianchi}
S.~Bubeck and N.~Cesa{-}Bianchi, ``Regret analysis of stochastic and
  nonstochastic multi-armed bandit problems,'' \emph{Foundations and Trends in
  Machine Learning}, vol.~5, no.~1, pp. 1--122, 2012.

\bibitem{gokcesu2018bandit}
K.~{Gokcesu} and S.~S. {Kozat}, ``An online minimax optimal algorithm for
  adversarial multiarmed bandit problem,'' \emph{IEEE Transactions on Neural
  Networks and Learning Systems}, vol.~29, no.~11, pp. 5565--5580, 2018.

\bibitem{tnnls3}
H.~R. Berenji and P.~Khedkar, ``Learning and tuning fuzzy logic controllers
  through reinforcements,'' \emph{IEEE Transactions on Neural Networks},
  vol.~3, no.~5, pp. 724--740, Sep 1992.

\bibitem{gokcesu2021generalized}
K.~Gokcesu and H.~Gokcesu, ``Generalized huber loss for robust learning and its
  efficient minimization for a robust statistics,'' \emph{arXiv preprint
  arXiv:2108.12627}, 2021.

\bibitem{ozkan}
H.~Ozkan, M.~A. Donmez, S.~Tunc, and S.~S. Kozat, ``A deterministic analysis of
  an online convex mixture of experts algorithm,'' \emph{IEEE Transactions on
  Neural Networks and Learning Systems}, vol.~26, no.~7, pp. 1575--1580, July
  2015.

\bibitem{gokcesu2020recursive}
K.~Gokcesu and H.~Gokcesu, ``Recursive experts: An efficient optimal mixture of
  learning systems in dynamic environments,'' \emph{arXiv preprint
  arXiv:2009.09249}, 2020.

\bibitem{singer}
A.~C. Singer and M.~Feder, ``Universal linear prediction by model order
  weighting,'' \emph{IEEE Transactions on Signal Processing}, vol.~47, no.~10,
  pp. 2685--2699, Oct 1999.

\bibitem{tnnls1}
R.~Song, F.~L. Lewis, and Q.~Wei, ``Off-policy integral reinforcement learning
  method to solve nonlinear continuous-time multiplayer nonzero-sum games,''
  \emph{IEEE Transactions on Neural Networks and Learning Systems}, vol.~PP,
  no.~99, pp. 1--10, 2016.

\bibitem{gokcesu2021optimally}
K.~Gokcesu and H.~Gokcesu, ``Optimally efficient sequential calibration of
  binary classifiers to minimize classification error,'' \emph{arXiv preprint
  arXiv:2108.08780}, 2021.

\bibitem{tnnls4}
J.~Moody and M.~Saffell, ``Learning to trade via direct reinforcement,''
  \emph{IEEE Transactions on Neural Networks}, vol.~12, no.~4, pp. 875--889,
  Jul 2001.

\bibitem{gokcesu2019outlier}
K.~Gokcesu, M.~M. Neyshabouri, H.~Gokcesu, and S.~S. Kozat, ``Sequential
  outlier detection based on incremental decision trees,'' \emph{{IEEE} Trans.
  Signal Process.}, vol.~67, no.~4, pp. 993--1005, 2019.

\bibitem{bartlett2019simple}
P.~L. Bartlett, V.~Gabillon, and M.~Valko, ``A simple parameter-free and
  adaptive approach to optimization under a minimal local smoothness
  assumption,'' in \emph{Algorithmic Learning Theory}.\hskip 1em plus 0.5em
  minus 0.4em\relax PMLR, 2019, pp. 184--206.

\bibitem{grill2015black}
J.-B. Grill, M.~Valko, and R.~Munos, ``Black-box optimization of noisy
  functions with unknown smoothness,'' \emph{Advances in Neural Information
  Processing Systems}, vol.~28, pp. 667--675, 2015.

\bibitem{piyavskii1972algorithm}
S.~Piyavskii, ``An algorithm for finding the absolute extremum of a function,''
  \emph{USSR Computational Mathematics and Mathematical Physics}, vol.~12,
  no.~4, pp. 57--67, 1972.

\bibitem{shubert1972sequential}
B.~O. Shubert, ``A sequential method seeking the global maximum of a
  function,'' \emph{SIAM Journal on Numerical Analysis}, vol.~9, no.~3, pp.
  379--388, 1972.

\bibitem{jacobsen1978global}
S.~E. Jacobsen and M.~Torabi, ``A global minimization algorithm for a class of
  one-dimensional functions,'' \emph{Journal of Mathematical Analysis and
  Applications}, vol.~62, no.~2, pp. 310--324, 1978.

\bibitem{mayne1984outer}
D.~Q. Mayne and E.~Polak, ``Outer approximation algorithm for nondifferentiable
  optimization problems,'' \emph{Journal of Optimization Theory and
  Applications}, vol.~42, no.~1, pp. 19--30, 1984.

\bibitem{mladineo1986algorithm}
R.~H. Mladineo, ``An algorithm for finding the global maximum of a multimodal,
  multivariate function,'' \emph{Mathematical Programming}, vol.~34, no.~2, pp.
  188--200, 1986.

\bibitem{shen1987interval}
Z.~Shen and Y.~Zhu, ``An interval version of shubert's iterative method for the
  localization of the global maximum,'' \emph{Computing}, vol.~38, no.~3, pp.
  275--280, 1987.

\bibitem{horst1987convergence}
R.~Horst and H.~Tuy, ``On the convergence of global methods in multiextremal
  optimization,'' \emph{Journal of Optimization Theory and Applications},
  vol.~54, no.~2, pp. 253--271, 1987.

\bibitem{breiman1993deterministic}
L.~Breiman and A.~Cutler, ``A deterministic algorithm for global
  optimization,'' \emph{Mathematical Programming}, vol.~58, no.~1, pp.
  179--199, 1993.

\bibitem{baritompa1994accelerations}
W.~Baritompa and A.~Cutler, ``Accelerations for global optimization covering
  methods using second derivatives,'' \emph{Journal of Global Optimization},
  vol.~4, no.~3, pp. 329--341, 1994.

\bibitem{hansen1995lipschitz}
P.~Hansen and B.~Jaumard, ``Lipschitz optimization,'' in \emph{Handbook of
  global optimization}.\hskip 1em plus 0.5em minus 0.4em\relax Springer, 1995,
  pp. 407--493.

\bibitem{sergeyev1998global}
Y.~D. Sergeyev, ``Global one-dimensional optimization using smooth auxiliary
  functions,'' \emph{Mathematical Programming}, vol.~81, no.~1, pp. 127--146,
  1998.

\bibitem{ellaia2012modified}
R.~Ellaia, M.~Z. Es-Sadek, and H.~Kasbioui, ``Modified piyavskii’s global
  one-dimensional optimization of a differentiable function,'' \emph{Applied
  Mathematics}, vol.~3, pp. 1306--1320, 2012.

\bibitem{brent2013algorithms}
R.~P. Brent, \emph{Algorithms for minimization without derivatives}.\hskip 1em
  plus 0.5em minus 0.4em\relax Courier Corporation, 2013.

\bibitem{gokcesu2021regret}
K.~Gokcesu and H.~Gokcesu, ``Regret analysis of global optimization in
  univariate functions with lipschitz derivatives,'' \emph{arXiv preprint
  arXiv:2108.10859}, 2021.

\bibitem{gokcesu2022low}
------, ``Low regret binary sampling method for efficient global optimization
  of univariate functions,'' \emph{arXiv preprint arXiv:2201.07164}, 2022.

\bibitem{horst2013global}
R.~Horst and H.~Tuy, \emph{Global optimization: Deterministic
  approaches}.\hskip 1em plus 0.5em minus 0.4em\relax Springer Science \&
  Business Media, 2013.

\bibitem{danilin1971estimation}
Y.~M. Danilin, ``Estimation of the efficiency of an absolute-minimum-finding
  algorithm,'' \emph{USSR Computational Mathematics and Mathematical Physics},
  vol.~11, no.~4, pp. 261--267, 1971.

\bibitem{malherbe2017global}
C.~Malherbe and N.~Vayatis, ``Global optimization of lipschitz functions,'' in
  \emph{International Conference on Machine Learning}.\hskip 1em plus 0.5em
  minus 0.4em\relax PMLR, 2017, pp. 2314--2323.

\bibitem{bouttier2020regret}
C.~Bouttier, T.~Cesari, and S.~Gerchinovitz, ``Regret analysis of the
  piyavskii-shubert algorithm for global lipschitz optimization,'' \emph{arXiv
  preprint arXiv:2002.02390}, 2020.

\end{thebibliography}

\end{document}